\documentclass[11pt]{amsart}
\usepackage{graphicx}
\usepackage{url}
\usepackage[round,sort]{natbib} % Use the natbib bibliography and citation package
\usepackage{vmargin}
\usepackage{breakurl}

\setmarginsrb{2cm}{2cm}{2cm}{2cm}{10mm}{10mm}{10mm}{10mm}
\theoremstyle{definition}

\theoremstyle{remark}

\numberwithin{equation}{section}
\begin{document}

\title[Calculus teaching and learning]{Calculus teaching and learning in South Korea}%
\author{Natanael Karjanto}%
\address{\hspace*{-0.55cm} Department of Mathematics \hfill Department of Mathematics \newline
University College \hfill School of Science and Technology \newline
Sungkyunkwan University \hfill Nazarbayev University \newline
Natural Science Campus, Suwon \hfill Astana \newline
South Korea \hfill Kazakhstan}%
\email{natanael@skku.edu \hfill natanael.karjanto@nu.edu.kz}%

\thanks{\textit{Jurnal Matematika Integratif} \textbf{9}(2): 179-193, 2013. Published by Department of Mathematics, Faculty of Mathematics and Natural Science, Padjadjaran University, West Java, Indonesia (ISSN 1412-6184).}%

%\date{}%
%\dedicatory{}%
%\commby{}%
% ----------------------------------------------------------------
\begin{abstract}
This article discusses an experience of teaching Calculus classes for the freshmen students enrolled at Sungkyunkwan University, one of the private universities in South Korea. The teaching and learning approach is a balance combination between the teacher-oriented traditional style of lecturing and other activities that encourage students for active learning and classroom participation. Based on the initial observation during several semesters, some anecdotal evidences show that students' learning is improved after implementing this student-oriented active learning approach, albeit a longer period of time is definitely needed to transform general students' attitude from passive learners to active ones.
\end{abstract}
\maketitle

% ----------------------------------------------------------------
\section{Introduction}
Many countries around the world are concerned with the quality of education of their citizenry, ranging from basic to higher education.
A number of international standards have been developed in
order to verify and compare the effects of education policies. Several well-known examples of assessments in
reading, mathematics, and science literacy amongst others are the Progress in International Reading Literacy
Study (PIRLS), the Trends in International Mathematics and Science Study (TIMSS),
and Programme for International Student Assessment (PISA). The first two studies were
conducted by the International Association for the Evaluation of Educational Achievement~\citep{Iea11}
in order to compare students' educational achievements across borders among
the participating nations. The latter one is coordinated by the Organization for Economic
Cooperation and Development~\citep{Oecd61} aiming at improving the educational policies and
outcomes not only for the OECD member countries but also for other non-OECD participating countries and economies.
The assessment results can help educational policy makers to shape their policy making in general and
to improve education strategies in particular. Although the assessments are predominantly targeted at the elementary and
secondary school levels, it is no doubt that the corresponding results may also reflect the quality of higher education.
In the realm of higher and college education, many institutions worldwide develop institutional or nationwide programmes
(or `visions') in teaching and learning to improve the quality of their higher education.

A lot of resources have been invested into these programmes to improve the quality of higher education, for instance
by training instructors, tailoring curricula, and introducing novel, non-traditional and innovative styles in
teaching and learning. Although the conventional-traditional teaching method still remains popular around
the globe, pioneer-minded institutions challenge these teaching methods by successfully
implementing up-to-date, non-traditional approaches in instruction, teaching and learning. Two well-known examples are `problem-based
learning'~\citep{Neville09,Hmelo04} and `guided-discovery'~\citep{Leutner93,Shrager86}, for which both are based on constructivist teaching
strategies. These two pedagogical approaches can be categorized as active learning activities, methods that provide harmony between the instructor's
role and the students' role.

This article describes an author's experience of teaching and learning Calculus in a four-year university in South Korea
and the discussion focuses on implementing a mixture teaching and learning methods between teacher-oriented and student-oriented approaches.
The content of the article is organized as follows. In the following section, the theoretical framework on active learning is laid and the literature is reviewed.
In Section~\ref{skku}, brief profiles of South Korea's higher education and Sungkyunkwan University (SKKU) are presented.
This university is a private institution in South Korea where the study has been conducted.
Next, Section~\ref{organ} explains the organization of Calculus courses at SKKU. This section describes the course management and the adopted textbook.
Section~\ref{ctl} follows with the current practice of Calculus teaching and learning at SKKU.
This section also discusses both the existing, traditional teaching method and the one that encourages active learning and student's participation. The final section draws conclusion to our discussion.

\section{Active learning}

There is an abundance source of literature on active learning, stimulating active participation,
and enhancing critical thinking and the list is not exhaustive. In the context of this article, active learning refers to techniques where students do more than simply listening to a lecture from an instructor.
The students are involved in doing some activities, to discover, to process, and to apply information~\citep{McKinney11}. The author offered several in-class examples of active learning techniques that can be used in both small and large classes, as well as with all level of students, amongst others are collaborative learning group, games, student debates, case study analysis, concept mapping, think-pair-share and student-led review sessions. Other concrete examples on active learning activities can also be found in~\citep{King93}.

\citet{Bonwell91} popularized the approach of active learning for class instruction.
They suggested that pupils can work in pairs, discuss the materials
while role-playing, debate, engage in case study, take part in cooperative learning, produce
short written exercises, and many other activities. \citet{Lazarus99} mentions a number of strategies
to stimulate active participation, for instance by questioning techniques and stimulating
group discussions. \citet{Kim06} explore students' participation in online
discussions. The authors discover that the participation patterns are
mainly influenced by the students' voluntary participation.
According to \citet{Meyers93}, active learning is derived from two basic assumptions: learning
by nature is an active endeavor, and different people learn in different ways. The authors
also show that the quantity of learning is enhanced when the students are engaged in active
learning. \citet{Salter10} emphasized that the implementation of problem-based learning techniques can change surface learning into
deep learning, memorized knowledge into understanding, etc.

Implementing active learning approach in teaching is not only beneficial for the students, but also useful for instructors, notably the new ones, as addressed by~\citet{Eison90}. In particular,
an attempt to encourage students who enrolled in Calculus classes to be active learners has been reported, amongst others by~\citet{Abramovitz10,Idris09} and \citet{Zhang03}.
It is interesting that even a Calculus textbook specializes on active learning edition has been published recently~\citep{Hughes-Hallett13}.
Active learning method has also been implemented not only in Calculus course, but also in other mathematics-related subjects, such as General Physics~\citep{Laws99} and Engineering Statistics~\citep{Kvam00}.
Related to the latter one, an evaluation of an active learning approach in teaching Introductory Statistics class has been reported by~\citet{Carlson11}.
The review of their article in the context of changed attitudes and improved performance is discussed recently in a newsletter edited by~\citet{Weiner12}.

As mentioned previously, the focus of this article is to share the author's experience in exerting active learning strategies toward the students enrolled in Calculus classes at Sungkyunkwan University (SKKU), a private, four-year university in South Korea.
Although this experience has led to some interesting anecdotal evidences, the author is fully convinced that this is not the end of of the journey toward improvements in teaching and learning.
It is hoped that after experiencing the active learning approach, the students will not only improve their learning skills in the mathematics related subjects, but also to elevate their learning achievements in other subjects along the road during their undergraduate study time frame and beyond.

In the following section, brief profiles of South Korea's higher education and SKKU are presented.

\section{Brief profiles of South Korea's higher education and Sungkyunkwan University} \label{skku}

\subsection{A brief profile of South Korea's higher education}
Education, including higher education, has a unique place in the South Korean's society. Thanks to its Confucian culture, the society puts a heavy emphasis on the value of gaining an education.
According to the World Factbook published by the Central Intelligence Agency, the literacy rate\footnote{The definition of literacy is a person of age 15 and over who can read and write.} in South Korea is 97.9\%~\citep{cia02}. Although South Korean students often rank extraordinarily in many international competitions, such as mathematics, science and informatics olympiads as well as comparative assessments, such as TIMSS, PIRLS and PISA, the education system in South Korea is heavily criticized since it put too much emphasis on passive learning and rote memorization~\citep{Chalmers97,Ho09}. This is one of the motivations to implement active teaching and learning in Calculus class and to conduct this research.

Although some form of higher education has existed continuously in South Korea since the fourth century, many universities have their history dated back prior the Japanese Imperial Period, such as Yonsei University Medical School \textsl{Gwanghyewon} (1885), \textsl{Ewha Hakdang} mission school for girls (1886) and Bosung College (1905), where the latter one is the embryo of Korea University.
After its liberation from the Japanese colonial rule in 1945, many of these universities are officially recognized and many more institutions are established: Ewha Womans University and Korea University (1945), Seoul National University, Yonsei University and SKKU (1946).

Nowadays, there are more than 370 official higher education institutions in the country, which includes more than 40 national universities, 180 private four-year universities and other types of higher education institutions, including polytechnics, cyber-universities, two- and three-year junior colleges~\citep{Parry11}. The advancement rate from secondary school to higher education in South Korea reaches almost 88\% in 2008, one among the highest in the OECD countries~\citep{Lee09}. The following subsection discusses a brief profile of SKKU.

\subsection{A brief profile of Sungkyunkwan University}
Sungkyunkwan University, often abbreviated as SKKU and also known by its nickname in Korean as \textsl{Seongdae}, is a private research university in the Republic of Korea. The university has two campus locations: the Humanities and Social Science Campus in Seoul, the country's capital and the Natural Science Campus in Suwon, which lies about 30~km south of Seoul and the provincial capital of Gyeonggi-do province. It boasts as the oldest university in the country, and possibly in the peninsula as the historical evidence shows that the chronicle of the university dated back as early as 1398 during the Joseon Dynasty. It is of no surprise that the former slogan of the university is `unique origin, unique future'.
The literal meaning of `Sungkyunkwan' itself is `an institution for building a harmonious society of perfected human beings'~\citep{Wiki}.
Fast forward to 1946, the current modern day university was established after the country's liberation form the Japanese Imperialism.
After half a century period, the Samsung Group provided funds to reacquire an almost depleted in resources university foundation in 1996~\citep{SKKU}.

\begin{center}
\begin{table}[h]
\begin{tabular}{|l|c|c|c|c|}
  \hline \hline
  % after \\: \hline or \cline{col1-col2} \cline{col3-col4} ...
      & 2013 & 2012        & 2011            & 2010   \\
  University  & International/ & International/ & International/ & International/   \\
  rankings    & National       & National       & National       & National   \\ \hline
  QS   & 162/6       & 179/6        & 259/7       & --     \\ \hline
  THE  & 201--225/6 & 201--225/5  & 301-350/6 & --    \\ \hline
  ARWU & 201--300/2--4    & 201--300/2--4    & 301--400/4--7   & 301--400/5--7  \\
  \hline \hline
\end{tabular}
  \vspace*{0.25cm}
  \caption{\small Sungkyunkwan University's rankings for the past couple of years according to the Quacquarelli Symonds (QS) World University ranking, Times of Higher Education (THE) and Academic Ranking of World Universities (ARWU). } \label{rank}
\end{table}
\end{center}
During the past several years, SKKU and other top universities in South Korea have been steadily showing their existence in the map of higher education worldwide.
According to Quacquarelli Symonds (QS) World University ranking, SKKU takes the 162$^\textmd{nd}$ position internationally and the 6$^\textmd{th}$ nationally in 2013.
According to Times Higher Education (THE), the university ranks in the category 201-225 for two consecutive years (2012 and 2013) and according to Academic Ranking of World Universities (ARWU), it ranks in the category of 201-300 for the past two years. ARWU is released by the Center for World-Class Universities at Shanghai Jiao Tong University. The university rankings from several sources for the past two to three years are listed in Table~\ref{rank}. Nowadays, SKKU has five colleges, ten divisions, 17 graduate schools and offers 51 majors~\citep{Kcue11}.
In addition, the university also boasts with 18,000 undergraduate students and 7000 graduate students, supported by more than 1000 faculty and 400 administrative staffs.

\section{Organization of Calculus course} \label{organ}

\subsection{Course management}
Two Calculus courses (Calculus~1 and Calculus~2) are offered as the Basic Science and Mathematics (BSM) modules in the Natural Science Campus at SKKU.
Other courses that belong BSM modules amongst others are General Physics 1 and 2, General Chemistry 1 and 2 and General Biology 1 and 2.
These BSM modules are offered and organized by the Faculty of Natural Sciences in cooperation with University College and the respective departments related to those courses.
A limited number of BSM modules are also offered in the Humanities and Social Science Campus in Seoul. Several faculty who reside in Suwon might need to commute to Seoul if they are assigned to teach in Seoul campus.

The two Calculus courses are by far the largest classes in terms
of student's number registered at SKKU. During Fall/Autumn 2010 semester, there were almost 1600
students enrolled in Calculus~2 courses and the majority of them are freshmen. A total of 18 sections are assigned to
12 instructors and one instructor acts as a course coordinator. Fourteen sections have a maximum number of students of 80 pupils and
four sections have larger size number of students up to a maximum of 120 pupils. All students would sit for identical examination
papers, both for the midterm and the final exams. However, each instructor may have additional
assessments, including but not limited to, homework, assignments, group projects or quizzes tailored to his/her own teaching.
The students may access the course syllabus through an intranet university's communication system called \textsl{icampus}~\citep{icampus}
as well as the internal registration system known as ASIS-GLS information system.

The semester at SKKU runs for 16 weeks, including a one-week period each for the midterm and the final exams. So, in principle there are only
14 weeks of effective teaching period during the entire semester. Each week teaching session is typically composed of either double one-and-half hour
sessions or a single three-hour session, depending on the timetable organized by the Mathematics Department.
During intensive semesters (Summer and Winter), the semester is organized in a three-week period, forming it to a total of 15 days which correspond to 15 weeks in a regular semester.

In particular, during Winter 2010 semester, the Department of Mathematics offers two sections of
Calculus 2 and two different instructors conduct both classes independently. As a
consequence, all exams are not necessarily identical due to this independent nature. Since the intensive semester is relatively short,
there are neither graded assignments nor group projects. There exist, however,
a few voluntary and ungraded assignments and group projects to gauge the students' interest
in learning the content knowledge. Similar to the regular semester, the maximum number of students in one class is
also limited to a maximum number of 80 pupils.

\subsection{Adopted textbook}
A single textbook is adopted and being used for teaching and learning of Calculus~1 and Calculus~2 courses, both in regular
(Spring and Fall/Autumn) as well as during intensive semesters (Winter and Summer). Although there is a plan to adopt
a different textbook for Calculus courses starting the new academic year 2011, so far only a particular textbook from~\citet{Stewart10} has been used for the two coursess.
The Korean translation of this very book is also readily available to the students and some of them prefer to read and to study the
translated version instead of the original one for a better understanding of the material. It is observed that many students made photocopies of
the book instead of purchasing the original version of it. It seems that violating of copyright law is still relatively rampant in this advanced country.

The textbook from~\citet{Stewart10} is intended for a year long course, covering both single and
multivariable Calculus, including an introduction to Vector Calculus. To obtain a general
overview of the syllabus for these two Calculus courses, simply divide the textbook
into two parts. Calculus~1 covers Chapter~1 to Chapter~8 and Calculus~2 covers Chapter~9
to Chapter~13 of the textbook. Although more chapters are weighted toward Calculus~1 course,
the pace of teaching can be done faster since the students are already knowledgeable already with
many elements of Calculus from their high school mathematics lessons~\citep{Kim13}. On the other hand, even though only
five chapters are covered in Calculus 2, there are occasions that the materials are covered relatively fast and
some parts are skipped entirely in order to explain better on other parts, which are rather new
to the students. In fact, the course coordinator has provided a detailed syllabus, including
the parts that need to be skipped and to be omitted.

Similar to many textbooks for basic courses at the freshman level of undergraduate programmes,
the adopted Calculus textbook~\citep{Stewart10} provides excellent materials for Calculus teaching and learning.
It presents a mixed balance of theoretical concepts, worked examples, use of technology and short biographies of several prominent mathematicians who contributed to the development of Calculus.
What is very essential for the students is its feature in providing many exercise problems
with diverse level of difficulties. The book is intended to be a mainstream Calculus text
that is suitable for every kind of course at every level. It is designed particularly for the
standard course of two semesters for students majoring in Natural Sciences and Engineering.
Certainly, students are expected to have some background of high school mathematics in
algebra, geometry and trigonometry as a general prerequisite for the course and in order to be able to follow the material in the book.
The book also contains many beautiful artworks that give a vivid illustration to the readers when particular concepts are being explained.
For every published edition, the author always makes an endeavour to improve the textbook, which includes
the mathematical precision, accuracy, clarity of exposition and outstanding examples and
problem sets which have also characterized the earlier editions too. The `focus on problem solving' feature
has made the book a favourite among students and instructors in a wide variety of colleges and
universities throughout the world. 

The applications of Calculus in many diverse fields of study are
also presented in an interesting manner. Since the book has a relatively challenging way of presentation, the students
who are rather weak in mathematics may have a difficult time following the textbook's rigorous approach and some
algebraic or conceptual jumps. The first three editions of the textbook were intended to be a synthesis of Calculus reform and
traditional approaches to Calculus instruction. In the fourth edition, the author continues to
follow that path by emphasizing conceptual understanding through visual, verbal, numerical
and algebraic approaches. The aim is to convey to the students both the practical power
of Calculus and the intrinsic beauty of the subject. In whatever teaching approaches that
one as an instructor may implement, or in whatever style approach the book might be presented, a
common goal remains to be achieved: to enable students to understand and appreciate the beauty of
Calculus. Some features of the book consists of conceptual exercises, graded exercise sets,
real-world data, applied projects, laboratory projects, writing projects, discovery projects,
problem solving, some rigourous proofs and the use of technology.

There are also ancillaries for both the instructors and the students. For instructors, the
PowerLecture CD-ROM is provided. This includes lectures in Power Point format and an
electronic version of the Instructor's Guide. Apart from this, the instructor also has an access
to the accompanying Complete Solutions Manual of the textbook and Tools for Enriching{\tiny\texttrademark}
Calculus available at \url{www.stewartcalculus.com}~\citep{Stewart05}. Ancillaries for students include Stewart
Specialty Website~\citep{Stewart05}, Enhanced WebAssign, The Brooks/Cole
Mathematics Resource Center Website~\citep{Cengage13}, Maple CD-ROM, Tools
for Enriching Calculus{\tiny \texttrademark}, Study Guide for both Single and Multivariable Calculus, Student
Solutions Manual and \textsl{Calc-Labs} with \textsl{Maple} and \textsl{Mathematica}.

\section{Calculus teaching and learning} \label{ctl}

\subsection{Teacher-oriented teaching and learning}
The teaching method implemented for Calculus courses, and any other mathematics based
courses, is generally traditional, conventional and teacher-oriented: the instructor acts as the controller of the
learning environment~\citep{Flinders13}. The instructor
of traditional way creates an environment for which the learning process is teacher-centered
instead of student-centered~\citep{Novak10}. Thus, in a number of occasions, the teaching
session is heavily dominated by merely delivering a lecture to the students with a minimal opportunity or even
no participation at all from the students' side. Nowadays, many experts in education around
the world encourage many educators, including school teachers and university instructors
to combine teaching methods, if not switching entirely.
Non-traditional pedagogical approaches such as problem-based learning and engaging students more by active learning and collaborative team work are definitely worth to consider.
There is an abundant literature on these topics, amongst others
are~\citep{Rosenthal95,Roj01,Zhang02,Ramsay06,Zwek06}.

The fact that many educators suggested to move away the teaching methods from the conventional ones, does
not necessarily mean that the traditional teaching method is not useful any longer. Depending on how the instructors
implement the traditional teaching method, an effective teaching that creates a conducive
learning environment and gives an optimal benefit for the students can still be achieved. By having a solid
preparation, implementing a structured organization, generating an outline, analyzing the audience, choosing
examples, choosing learning activities, using audio/visual aids and reviewing the materials,
interesting set of lectures can be created using this traditional teaching approach~\citep{Rutgers13}.
Although some might believe that non-conventional
or modern teaching methods are better than the traditional and conventional ones, alternative
methods are equally successful when they are handled in an effective manner. 
In Calculus courses, entirely abandoning the traditional and conventional teaching method can be very challenging.
However, elements of non-conventional teaching methods can be embedded so
that the teaching activities do not only have more variations but also become more interesting
for the students. In turn, it is expected that the students will have a voracious appetite
in learning and studying.

The PowerLecture CD-ROM from the textbook is fully utilized during the teaching activities. When explaining
new concepts, definitions, properties and theorems, the CD-ROM is very useful since it saves
time compared to re-writing again everything on the blackboard/whiteboard. The CD-ROM
also contains good quality artworks that describe mathematical concepts or three-dimensional
objects which are not so easy or even impossible to sketch by hand. When discussing
examples, the CD-ROM provides a step by step explanation, identical to one found in the textbook.
The students can easily lose their attention during lecture, particularly when the lecture slide has been shown for more than 20 minutes.
So, instead of directly displaying the step by step explanation from the slide, it is good to explain the students by hand, either by writing
on the computer screen (projected to a huge screen) or on the whiteboard. By doing this, hopefully we are able to keep the students'
attention and to improve their understanding in problem solving.
After this step, the students can be inquired whether they have any questions or find any unclear explanation steps.
It is a norm culturally that Korean students will keep quiet and will not pose any questions in public, even though they might not fully grasp the entire teaching session.
Assigning the students exercise problems from the textbook can maintain a high level of alertness from the students' side.

Every lecture session would be recorded (the so-called \textsl{e+ lecture}), and the recording is
uploaded automatically to the \textsl{icampus} system after the teaching session is completed.
The recording software does not only capture the explanation from the instructor but also includes the
slides and everything written on the computer screen. This is beneficial to the students since they
might want to review the lecture again after the class is over in case they find certain concepts are difficult to grasp and
the instructor's pace is too fast to follow and yet they are rather shy to ask questions in a public setting or to inquire
the instructor to repeat the explanation on certain topics, particularly difficult concepts. The recording is also
useful for the students who miss the class during that particular session of the day.

The classrooms used for Calculus 2 classes during Fall/Autumn and Winter 2010 are 31255 and
26515, respectively. Both classrooms are lecture hall type and the desk and chairs are
fixed. Unfortunately re-arranging seat, for the sake of some active learning activities, is impossible. Since the
students always sit at the back of the classroom or at the far side seats,
taking the front and the middle seats is always encouraged. Certainly, it is not ethical to enforce where the students must sit, but it is
preferable that as many students as possible use these seats and leave the rest for the latecomers.
The response to this admonition is rather poor. Although some students move to
the front seats, the majority still remained seated at the rear part of the classroom. Very often, the first two
or three rows of the front seats remain empty for the lecture hall 31255 and the first front row remains empty for
the lecture hall 26515. 

Furthermore, the majority of the student prefers to sit close to their peers whom they already familiar.
When the students sign up for a particular course, they could have some discussions and agreements
with their peers before selecting a particular instructor. This may explain why some students
prefer to mingle with their peers rather than to mix around with other `non-familiars'. On the
other hand, as the semester evolves, some students were making new friends with others whom
they did not really know initially. Signing up the class is on the first come first served based.
As a consequence, the composition of each class may vary not only on their major study but
also according to their level year of study.

\subsection{Student-oriented teaching and learning}
The student-oriented teaching and learning in Calculus classes has been implemented by active learning strategies.
Some of the adopted strategies are question and answer session, peer-tutoring session and team-work project assignment.
The students' active participation can be stimulated by posing them several questions during
the teaching session and allow them to give responses. For instance, it may happen that some calculation and typographical errors occur when discussing examples taken from the textbook.
The students may or may not be aware of the errors displayed on the whiteboard or the projected screen. Even those who are aware, they are hesitant to pinpoint the errors publicly.
This situation is a good opportunity to pause and ask them to pinpoint where the errors are.
It helps them to think, to stay alert and to be involved during the classroom session.

Another way is by asking one or more students to write and present their work on the whiteboard.
From this presentation, the way they present mathematical expressions can be monitored and suggestions for improvement are given.
The way they present on the whiteboard generally also reflects the way they write their solutions on the exams.
Even if it is no problem for them to verify the correct answers of the odd-numbered questions from the textbook, the presentation of the solutions
may need to be improved. Furthermore, many of the students find it difficult to explain their presentation in English even though
they can explain comfortably in Korean. So, after the students complete their work on the whiteboard, a further
training could be provided by inquiring them some questions of what they have written and let them explain
with their own words. It is often discovered that they write certain symbols and yet did not
explain what those symbols mean. This training is hopefully useful for them when eventually
they have to write their exam papers.

\section{Conclusion}

\subsection{Weakness}

Although embracing students' participation strategies in the teaching style is an indispensable
pedagogy tool, the students' psychology can easily throw a spanner in the works. If the
instructor does not specify well-defined regulations for the students' rewards for participating
in the activities, students would lose interest and their efforts would diminish as well. The
regulations must clearly state what grade incentives students can earn for their participation.
Students therefore need to know how this is evaluated and how it affects their grade.
Especially students who perform below average seem to need this kind of incentives.

\subsection{Anecdotal evidence}

From this study we have seen some anecdotal evidences that students' learning is improved when they are actively
engaged with the study material, instead of only sitting passively in the classroom and listening to the lecture. More
success can be achieved when the classroom activities are also fun. With a proper balance
between lecturing and engaging students new concepts and activities in which students, alone or in groups,
need to struggle themselves with these concepts, makes the learning time in the classroom
more effective and the time spent in class becomes more enjoyable for the students. Students
also show appreciation for this style of teaching. Despite classroom participation and group
activities are generally more successful in classes with a small number of students,
some promising and good results have been accomplished even though the class sizes were relatively large.

\subsection{Future study}

Improving students' participation in Calculus classes has achieved a measure of success to
a certain extent. Some students are willing to participate in class activities, but others are still quite hesitant.
Because the students' educational backgrounds are diverse, particularly the weaker students should discern how to benefit from classroom participation. Although it is almost impossible to design a teaching style that suits for all students, some student-oriented active learning approaches can be implemented for the future study, including team-work guided discovery~\citep{Leutner93,Shrager86,Dumitrascu09} or problem-based learning approach~\citep{Hmelo04,Neville09}.
Implementing practical methods that worked out best for the students will be implemented continuously and improved regularly during the coming semesters.

\section*{Acknowledgement}
The author acknowledged many fruitful discussions with Rob Lahaye from Physics Department, University College, Natural Science Campus of SKKU and Dan Kim from Center of Teaching and Learning, Humanities and Social Science Campus of SKKU and Educational Broadcasting System (EBS) television as well as
SKKU's Academic Affairs Office for a generous financial support through the Undergraduate Educational Policy Development Project.

% ----------------------------------------------------------------
%\bibliographystyle{amsplain}
%\bibliography{}

%\bibliographystyle{apacite}

\end{document}